\documentclass[12pt]{article}
\usepackage{amssymb}
\usepackage{amsmath}

\begin{document}

\begin{center}
\textbf{\large The Sobolev space $W_2^{1/2}$: Simultaneous improvement of 
functions by a homeomorphism of the circle} 
\end{center}

\begin{center}
\textsc{Vladimir Lebedev}
\end{center}

\quad

\begin{quotation}
{\small \textbf{Abstract.} It is known that for every continuous real-valued 
function $f$ on the circle $\mathbb T=\mathbb R/2\pi\mathbb Z$ there exists a 
change of variable, i.e., a self-homeomorphism $h$ of $\mathbb T$, such that 
the superposition $f\circ h$ is in the Sobolev space $W_2^{1/2}(\mathbb T)$. 
We obtain new results on simultaneous improvement of functions by a single 
change of variable in relation to the space $W_2^{1/2}(\mathbb T)$.  The main 
result is as follows: there does not exist a self-homeomorphism $h$ of 
$\mathbb T$ such that $f\circ h\in W_2^{1/2}(\mathbb T)$ for every $f\in 
\mathrm{Lip}_{1/2}(\mathbb T)$. Here $\mathrm{Lip}_{1/2}(\mathbb T)$ is the 
class of all functions on $\mathbb T$ satisfying the Lipschitz condition of 
order $1/2$. 


} 
\end{quotation}

\textbf{1. Introduction.}  Given an integrable function $f$ on the circle 
$\mathbb T=\mathbb R/2\pi\mathbb Z$, consider its Fourier expansion: 
$$
f(t)\sim\sum_{k\in\mathbb Z}\widehat{f}(k)e^{ikt}, \qquad t\in\mathbb T.
$$
Recall that the Sobolev space $W_2^{1/2}(\mathbb T)$ is the space of all 
square integrable functions $f$ satisfying 
$$
\|f\|_{W_2^{1/2}(\mathbb T)}=\bigg(\sum_{k\in\mathbb Z}|\widehat{f}(k)|^2|k|\bigg)^{1/2}<\infty.
\eqno(1)
$$
In what follows, by $C(\mathbb T)$ we denote the Banach space of all 
continuous complex-valued functions $f$ on $\mathbb T$ with the usual norm 
$\|\cdot\|_{C(\mathbb T)}$ defined by $\|f\|_{C(\mathbb T)}=\sup_{t\in\mathbb 
T}|f(t)|$. Given a modulus of continuity $\omega$, i.e., a nondecreasing 
continuous function on $[0, +\infty)$ with $\omega(0)=0$, by 
$\mathrm{Lip}_\omega(\mathbb T)$ we denote the class of all functions $f\in 
C(\mathbb T)$ satisfying $\omega(f, \delta)=O(\omega(\delta)), 
~\delta\rightarrow+0,$ where 
$$
\omega(f, \delta)=\sup_{|t_1-t_2|\leq\delta} |f(t_1)-f(t_2)|,
\qquad \delta\geq 0,
$$
is the modulus of continuity of $f$. For $0<\alpha\leq 1$ we just write 
$\mathrm{Lip}_\alpha$ instead of $\mathrm{Lip}_{\delta^\alpha}$. By 
$V(\mathbb T)$ we denote the space of all functions of bounded variation on 
$\mathbb T$. 

It is well known that one can improve certain properties of a continuous 
function related to its Fourier series by an appropriate change of variable, 
i.e., by a self-homeomorphism of $\mathbb T$. The first result on this 
subject is the Bohr--P\'al theorem stating that for every real-valued $f$ in 
$C(\mathbb T)$ there exists a self-homeomorphism $h$ of $\mathbb T$ such that 
the superposition $f\circ h$ belongs to the space $U(\mathbb T)$ of functions 
with uniformly convergent Fourier series. In fact, the original proof of Bohr 
and P\'al yields 
$$
f\circ h=g+\psi, 
\eqno(2)
$$ 
where $g\in W_2^{1/2}\cap C(\mathbb T)$ and $\psi\in V\cap C(\mathbb T)$. 
This implies $f\circ h\in U(\mathbb T)$, since $W_2^{1/2}\cap C(\mathbb T)$ 
and $V\cap C(\mathbb T)$ are subsets of $U(\mathbb T)$. 

Subsequently, for certain function spaces the question of whether every 
continuous function can be transformed by a suitable homeomorphic change of 
variable into a function that belongs to a given space was studied by various 
authors. Some of these studies concern the possibility of simultaneous 
improvement of several functions by means of a single change of variable.  
For a survey on the subject see [3], [11]; later results were obtained in 
[2], [5--9].\footnote{In [5] $C(\mathbb T)$ stands for the space of all 
\emph{real-valued} continuous functions on $\mathbb T$.} It is worth noting 
that the proof by Bohr and P\'al is based on the Riemann's theorem on 
conformal mappings, while investigations that followed mostly involve real 
analysis methods. 

We note, in particular, that an immediate consequence of [12, Corollary 1] 
(see also [5, Sec. 3], [2], [8]) is that the term $\psi$ in (2) can be 
omitted, i.e., the following refined version of the Bohr--P\'al theorem 
holds: for every real-valued $f\in C(\mathbb T)$ there exists a 
self-homeomorphism $h$ of $\mathbb T$ such that $f\circ h\in 
W_2^{1/2}(\mathbb T)$.  

The first result on simultaneous improvement was obtained in [4] in relation 
to the space $U(\mathbb T)$: if $K$ is a compact set in $C(\mathbb T)$, then 
there exists a self-homeomorphism $h$ of $\mathbb T$ such that $f\circ h\in 
U(\mathbb T)$ for every $f\in K$. In other words, given a modulus of 
continuity $\omega$, there is a change of variable $h$ such that $f\circ h\in 
U(\mathbb T)$ for every $f\in \mathrm{Lip}_\omega(\mathbb T)$. This result 
naturally led to a question if it is possible to attain the condition $f\circ 
h\in W_2^{1/2}(\mathbb T)$ for every $f\in K$. A negative answer was obtained 
in [5, Theorem 4]; as it turned out, given a real-valued $u\in C(\mathbb T)$, 
the property that for every real-valued $v\in C(\mathbb T)$ there is a 
homeomorphism $h$ such that both $u\circ h$ and $v\circ h$ are in 
$W_2^{1/2}(\mathbb T)$ is equivalent to the boundness of variation of $u$. 
Thus, in general, there is no single change of variable which will bring two 
real-valued functions in $C(\mathbb T)$ into $W_2^{1/2}(\mathbb T)$. This 
amounts to the existence of a complex-valued $f\in C(\mathbb T)$ such that 
$f\circ h\notin W_2^{1/2}(\mathbb T)$ whenever $h$ is a self-homeomorphism of 
$\mathbb T$. 

In the present paper we obtain new results on simultaneous improvement of 
functions by a single change of variable in relation to the space 
$W_2^{1/2}(\mathbb T)$. 

For an arbitrary function $F\in C(\mathbb T)$ consider the family 
$S_F=\{F(\cdot\,+\theta) : \theta\in\mathbb T\}$ of translations of $F$ and 
the family $K_F=\{F\ast\lambda : \lambda\in P(\mathbb T)\}$  of convolutions 
of $F$ with probability measures ($P(\mathbb T)$ is the set of all 
probability measures on $\mathbb T$). Both  $S_F$ and $K_F$ are compact sets 
in $C(\mathbb T)$ and $S_F\subseteq K_F$. Clearly, if $\lambda\in P(\mathbb 
T)$, then $|\widehat{F\ast\lambda}(k)|\leq |\widehat{F}(k)|, \,k\in\mathbb 
Z$. Thus, if $F$ is in $W_2^{1/2}(\mathbb T)$, then $S_F$ and $K_F$ are just 
subsets of  $W_2^{1/2}(\mathbb T)$. In Section 2 we show that if $F$ is not 
in $W_2^{1/2}(\mathbb T)$, then, under certain additional assumption on $F$, 
there does not exist a homeomorphism which will bring all functions in $S_F$ 
into $W_2^{1/2}(\mathbb T)$. In Section 3, without any additional assumptions 
on $F$, we obtain a similar result for $K_F$. 

In Section 4 we obtain the main result of the paper: there does not exist a 
self-homeomorphism $h$ of $\mathbb T$ such that $f\circ h\in 
W_2^{1/2}(\mathbb T)$ for every $f\in\mathrm{Lip}_{1/2}(\mathbb T)$. This is 
a direct consequence of the result on translations (as well as of the result 
on convolutions) and the well-known fact that $\mathrm{Lip}_{1/2}(\mathbb 
T)\nsubseteq W_2^{1/2}(\mathbb T)$. We note that in [7] it was shown that if 
$\alpha<1/2$, then there exist two real-valued functions in 
$\mathrm{Lip}_\alpha(\mathbb T)$ such that there is no single change of 
variable which will bring them into $W_2^{1/2}(\mathbb T)$. The author does 
not know if such a pair of functions can be found in 
$\mathrm{Lip}_{1/2}(\mathbb T)$. It is also worth noting that for 
$\alpha>1/2$ functions in $\mathrm{Lip}_\alpha(\mathbb T)$ do not require an 
improvement, since in this case $\mathrm{Lip}_\alpha(\mathbb T)\subseteq 
W_2^{1/2}(\mathbb T)$. The imbedding follows from the well-known equivalence 
of the seminorm  $\|\cdot\|_{W_2^{1/2}(\mathbb T)}$ (see (1)) and the 
seminorm $\||\cdot|\|_{W_2^{1/2}(\mathbb T)}$ defined by 
$$
\||f|\|_{W_2^{1/2}(\mathbb T)}=
\bigg(\int_{[0, 2\pi]}\frac{1}{t^2}\bigg(\int_{[0, 2\pi]}|f(x+t)-f(x)|^2 dx\bigg) dt\bigg)^{1/2}.
\eqno(3)
$$
 
The concluding Section 5 contains certain remarks, open problems and the 
shortest, known to the author, proof of the refined version of the 
Bohr--P\'al theorem. 

\quad

\textbf{2. Translations of a continuous function of analytic type.} Let $F\in 
C(\mathbb T)$. For each $\theta\in\mathbb T$ define the function $F_\theta$ 
by $F_\theta(t)=F(t+\theta)$. Consider the family $S_F$ of translations of 
$F$: 
$$
S_F=\{F_\theta : \theta\in\mathbb T\}.
$$
(Clearly, $S_F$ is a compact set in $C(\mathbb T)$.) 

By $C^+(\mathbb T)$ we denote the class of all continuous functions of 
analytic type on $\mathbb T$, i.e., of those $F\in C(\mathbb T)$ which 
satisfy $\widehat{F}(k)=0$ for all $k<0$.

It is obvious that if $F\in W_2^{1/2}(\mathbb T)$, then $S_F\subseteq 
W_2^{1/2}(\mathbb T)$. On the other hand, the following theorem holds. 

\quad

\textbf{Theorem 1.} \emph{Let $F\in C^+(\mathbb T)$. Suppose that $F\notin 
W_2^{1/2}(\mathbb T)$. Then there does not exist a self-homeomorphism $h$ of 
$\mathbb T$ such that $f\circ h\in W_2^{1/2}(\mathbb T)$ for every $f\in 
S_F$.} 

\quad 

To prove the theorem we will need Lemma 1 below. It has a technical character 
and will also be used in the next section. Before we proceed to the lemma, 
note that the bilinear form 
$$
B(x, y)=\frac{1}{2\pi}\int_{\mathbb T}x(t)dy(t),
$$ 
defined for $x\in C(\mathbb T), \,y\in V(\mathbb T)$, is invariant (up to a 
sign) with respect to self-homeomorphisms $h$ of $\mathbb T$, namely, if 
$x\in C(\mathbb T), \,y\in V(\mathbb T)$, then $x\circ h\in C(\mathbb T),\ 
y\circ h\in V(\mathbb T),$ and $B(x\circ h, y\circ h)=\pm B(x, y)$, where the 
choice of a sign depends on whether $h$ preserves the orientation or reverses 
it. Note also that if $x$ and $y$ are in $W_2^{1/2}(\mathbb T)$, then 
$\sum_{k\in\mathbb Z}|\widehat{x}(-k)\,ik\, \widehat{y}(k)|<\infty$.  

\quad

\textbf{Lemma 1.} \emph{Let $x\in C(\mathbb T), \,y\in V(\mathbb T)$ and $x, 
y\in W_2^{1/2}(\mathbb T)$. Then} 
$$ 
\textrm{(i)}\qquad\qquad\qquad \frac{1}{2\pi}\int_{\mathbb T}x(t)dy(t)=\sum_{k\in\mathbb Z}\widehat{x}(-k)\,ik\, 
\widehat{y}(k); 
$$ 
$$ 
\textrm{(ii)}\qquad\qquad \bigg|\frac{1}{2\pi}\int_{\mathbb T}x(t)dy(t)\bigg|\leq 
\|x\|_{W_2^{1/2}(\mathbb T)}\|y\|_{W_2^{1/2}(\mathbb T)}. 
$$

\quad

\textbf{Proof.} Clearly, (ii) follows from (i). To verify (i), observe that 
$$
\frac{1}{2\pi}\int_{\mathbb T} e^{-ikt}dy(t)=-\frac{1}{2\pi}\int_{\mathbb T}y(t)d e^{-ikt}=
ik\widehat{y}(k),\qquad k\in\mathbb Z,
$$ 
so (i) holds in the case when $x$ is a trigonometric polynomial. In the 
general case it suffices to approximate $x$ by the Fej\'er sums 
$\sigma_N(x)$: 
$$
\sigma_N(x)(t)=\sum_{|k|\leq N}\widehat{x}(k)\bigg(1-\frac{|k|}{N}\bigg)e^{ikt}, \quad N=1, 2, \ldots.
\eqno(4)
$$ 
Indeed, we have
$$
\frac{1}{2\pi}\int_{\mathbb T}\sigma_N(x)(t)dy(t)=
\sum_{|k|\leq N} \widehat{x}(-k)\bigg(1-\frac{|k|}{N}\bigg)\,ik\,\widehat{y}(k).
\eqno(5)
$$
Let $N\rightarrow\infty$. The sequence $\sigma_N(x), \,N=1, 2, \ldots,$ 
converges uniformly to $x$, so
$$
\frac{1}{2\pi}\int_{\mathbb T}\sigma_N(x)(t)dy(t)\rightarrow \frac{1}{2\pi}\int_{\mathbb T} x(t)dy(t)
$$ 
and it remains to notice that the right-hand side in (5) tends to the 
right-hand side in (i). The lemma is proved. 

\quad

\textbf{Proof of Theorem 1.} We have $F\in C(\mathbb T)$ and
$$
F(t)\sim\sum_{n\geq 0} c_n e^{int},
$$
where, by the assumption,
$$
\sum_{n\geq 0} |c_n|^2 n=\infty.
\eqno(6)
$$

Suppose that, contrary to the assertion of the theorem, there exists a 
self-homeomorphism $h$ of $\mathbb T$ such that $F_\theta\circ h\in 
W_2^{1/2}(\mathbb T)$ for every $\theta\in\mathbb T$. Consider the sets 
$T_m\subseteq\mathbb T, \,m=1, 2, \ldots,$ defined by
$$
T_m=\{\theta\in\mathbb T : \|F_\theta\circ h\|_{W_2^{1/2}(\mathbb T)}\leq m\}.
$$
One can easily see that the sets $T_m, \,m=1, 2, \ldots,$ are closed. Indeed, 
given an $m$, assume that a sequence $\theta_n, n=1, 2, \ldots,$ is in $T_m$ 
and converges to a $\theta_0$. Then, for an arbitrary positive integer $p$ we 
have $\sum_{|k|\leq p}|\widehat{F_{\theta_n}\circ h}(k)|^2|k|\leq m^2$. Since 
the sequence $F_{\theta_n}\circ h, n=1, 2, \ldots,$ converges uniformly to 
$F_{\theta_0}\circ h$, we obtain $\sum_{|k|\leq p}|\widehat{F_{\theta_0}\circ 
h}(k)|^2|k|\leq m^2$. Since $p$ was chosen arbitrarily this implies 
$\theta_0\in T_m$. Hence, each $T_m$ is closed. At the same time 
$$
\mathbb T=\bigcup_{m=1}^\infty T_m,
$$
whence, using the Baire category theorem, we see that at least one of the 
sets $T_m$, say $T_{m_0}$, contains an interval $I\subseteq\mathbb T$. So 
$\|F_\theta\circ h\|_{W_2^{1/2}(\mathbb T)}\leq m_0$ for all $\theta\in I$. 
Replacing, if necessarily, $h$ with $h+\gamma_I$, where $\gamma_I$ is the 
center of $I$, one can assume that $I=(-\delta_0, \delta_0)$, where 
$0<\delta_0\leq\pi$. Thus, 
$$
\|F_\theta\circ h\|_{W_2^{1/2}(\mathbb T)}\leq m_0 \quad \textrm{for all} \quad \theta\in (-\delta_0, \delta_0).
\eqno(7)
$$

For $0<\delta<\delta_0$ we let
$$
F^\delta(t)=\frac{1}{2\delta}\int_{-\delta}^\delta F(t+\theta)d\theta.
$$
Note that for $x, t\in \mathbb T$ we have 
$$
|F^\delta\circ h(x+t)-F^\delta\circ h(x)|\leq
$$
$$
\leq\frac{1}{2\delta}\int_{-\delta}^\delta |F_\theta\circ h(x+t)-F_\theta\circ h(x)|d\theta\leq
\bigg(\frac{1}{2\delta}\int_{-\delta}^\delta |F_\theta \circ h(x+t)-F_\theta\circ h(x)|^2d\theta\bigg)^{1/2},
$$
so for an arbitrary $\varepsilon$ such that $0<\varepsilon<2\pi$, we obtain 
$$
\int_\varepsilon^{2\pi}\frac{1}{t^2}\bigg(\int_0^{2\pi}|F^\delta\circ h(x+t)-F^\delta\circ h(x)|^2dx\bigg)dt\leq
$$
$$
\leq\frac{1}{2\delta}\int_{-\delta}^\delta\bigg(\int_\varepsilon^{2\pi}\frac{1}{t^2}\bigg(\int_0^{2\pi}
|F_\theta\circ h(x+t)-F_\theta\circ h(x)|^2dx\bigg)dt\bigg)d\theta,
$$
whence, taking into account (7) and the equivalence of the seminorms\\ 
$\|\cdot\|_{W_2^{1/2}(\mathbb T)}$ and $\||\cdot|\|_{W_2^{1/2}(\mathbb T)}$ 
(see (3)), we obtain that $F^\delta\circ h\in W_2^{1/2}(\mathbb T)$ and 
$$
\|F^\delta\circ h\|_{W_2^{1/2}(\mathbb T)}\leq c \quad \textrm{for all} \quad \delta\in(0, \delta_0),
\eqno(8)
$$
where $c>0$ does not depend on $\delta$. 

It is clear that $F^\delta$ is continuous and of bounded variation. So 
$F^\delta\circ h$ is also continuous and of bounded variation. Using Lemma 1 
(see (ii)) and (8), we have that for all $\delta\in (0, \delta_0)$ 
$$
\bigg|\frac{1}{2\pi}\int_\mathbb T \overline{F^\delta(t)}dF^\delta(t)\bigg|=
\bigg|\frac{1}{2\pi}\int_\mathbb T \overline{F^\delta\circ h(t)}d(F^\delta\circ h)(t)\bigg|
\leq
$$
$$
\leq\|\overline{F^\delta\circ h}\|_{W_2^{1/2}(\mathbb T)}\|F^\delta\circ h\|_{W_2^{1/2}(\mathbb T)}\leq c^2,
\eqno(9)
$$
where the bar stands for the complex conjugation.

Let 
$$
\chi_\delta=\frac{1}{2\delta} 1_{(-\delta, \,\delta)},
$$ 
where $1_{(-\delta, \delta)}$ is the indicator function of the interval 
$(-\delta, \delta)$. Obviously, $F^\delta=F\ast\chi_\delta$ and 
$$
\widehat{\chi_\delta}(k)=\frac{1}{2\pi}\frac{\sin k\delta}{k\delta}, \qquad k\neq 0.
$$ 
Applying Lemma 1 (see (i)), we obtain 
$$
\frac{1}{2\pi}\int_\mathbb T \overline{F^\delta(t)}dF^\delta(t)=
\sum_{k\in\mathbb Z}\overline{\widehat{F^\delta}(k)}ik\widehat{F^\delta}(k)=
i\sum_{k\in\mathbb Z}|\widehat{F^\delta}(k)|^2k=
$$
$$
=i\sum_{k\in\mathbb Z}|\widehat{F}(k)|^2|\widehat{\chi_\delta}(k)|^2k
=i\sum_{n\geq 1}|c_n|^2\bigg(\frac{1}{2\pi}\frac{\sin n\delta}{n\delta}\bigg)^2n.
$$
Thus (see (9)), for all $\delta\in (0, \delta_0)$
$$
\sum_{n\geq 1}|c_n|^2\bigg(\frac{1}{2\pi}\frac{\sin n\delta}{n\delta}\bigg)^2n\leq c^2.
$$

Choose a positive integer $N$. We see that for all $\delta\in (0, \delta_0)$ 
$$
\sum_{n=1}^N|c_n|^2\bigg(\frac{\sin n\delta}{n\delta}\bigg)^2n\leq (2\pi c)^2.
$$
Let $\delta\rightarrow+0$. We obtain
$$
\sum_{n=1}^N|c_n|^2n\leq (2\pi c)^2,
$$
which contradicts (6), since $N$ was chosen arbitrarily. The theorem is 
proved. 

\quad

\textbf{3. Convolutions of a continuous function with probability measures.} 
Let $F\in C(\mathbb T)$. Let $P(\mathbb T)$ be the set of all probability 
measures on $\mathbb T$. Consider the family $K_F$ of convolutions of $F$ 
with the measures in $P(\mathbb T)$: 
$$
K_F=\{F\ast\lambda : \lambda\in P(\mathbb T)\}.
$$
Clearly, $K_F$ is a compact set in $C(\mathbb T)$. It is also clear that (as 
we mentioned in the introduction) if $F\in W_2^{1/2}(\mathbb T)$, then 
$K_F\subseteq W_2^{1/2}(\mathbb T)$. On the other hand, the following theorem 
holds. 

\quad

\textbf{Theorem 2.} \emph{Let $F\in C(\mathbb T)$. Suppose that $F\notin 
W_2^{1/2}(\mathbb T)$. Then there does not exist a self-homeomorphism $h$ of 
$\mathbb T$ such that $f\circ h\in W_2^{1/2}(\mathbb T)$ for every $f\in 
K_F$.} 

\quad 

To prove Theorem 2 we will need two lemmas. As in the previous section, for 
an $F\in C(\mathbb T)$ and a $\theta\in\mathbb T$, we use $F_\theta$ to 
denote the function defined by $F_\theta(t)=F(t+\theta)$. 

\quad

\quad

\textbf{Lemma 2.} \emph{Let $x\in C(\mathbb T)$. Let $\varphi$ be a 
continuous mapping of $\mathbb T$ into itself. Then for each $\nu\in\mathbb 
Z$ the function $\theta\rightarrow\widehat{x_\theta\circ \varphi}(\nu)$ is 
continuous on $\mathbb T$ and 
$$
\frac{1}{2\pi}\int_\mathbb T |\widehat{x_\theta\circ \varphi}(\nu)|^2
d\theta=\sum_{k\in\mathbb Z}|\widehat{x}(k)|^2|\widehat{e^{ik\varphi}}(\nu)|^2.
\eqno(10)
$$} 

\quad

\textbf{Proof.} The continuity of $\widehat{x_\theta\circ \varphi}(\nu)$ is 
obvious. To prove (10), assume first that $x$ is a trigonometric polynomial. 
Then (there is a finite number of nonzero terms in the sum below) 
$$
x_\theta(t)=x(t+\theta)=\sum_{k\in\mathbb Z}\widehat{x}(k)e^{ik(t+\theta)},
$$
whence
$$
x_\theta\circ \varphi(t)=\sum_{k\in\mathbb Z}\widehat{x}(k)e^{ik\varphi(t)}e^{ik\theta}.
$$
So
$$
\widehat{x_\theta\circ \varphi}(\nu)=\sum_{k\in\mathbb Z}\widehat{x}(k)\widehat{e^{ik\varphi}}(\nu)e^{ik\theta},
$$
which implies (10). In the general case, consider the Fej\'er sums 
$\sigma_N(x), N=1, 2, \ldots,$ (see (4)). We have 
$$
\frac{1}{2\pi}\int_\mathbb T |(\sigma_N(x)_\theta\circ \varphi)^\wedge(\nu)|^2
d\theta=\sum_{|k|\leq N}|\widehat{x}(k)|^2\bigg(1-\frac{|k|}{N}\bigg)^2|\widehat{e^{ik\varphi}}(\nu)|^2.
$$ 
It remains to tend  $N$ to $\infty$. The lemma is proved. 

\quad 

\textbf{Lemma 3.} \emph{Let $h$ be a self-homeomorphism of $\mathbb T$. Let
$k\in\mathbb Z$. Suppose that $e^{ikh}\in W_2^{1/2}(\mathbb T)$. Then
$\|e^{ikh}\|_{W_2^{1/2}(\mathbb T)}\geq |k|^{1/2}$.} 

\quad 

\textbf{Proof.} Applying Lemma 1 (see (ii)), we obtain
$$
|k|=\bigg|\frac{1}{2\pi}\int_\mathbb T e^{ik x} de^{-ik x}\bigg|=
\bigg|\frac{1}{2\pi}\int_\mathbb T e^{ik h(t)} de^{-ik h(t)}\bigg|\leq \|e^{ikh}\|_{W_2^{1/2}(\mathbb T)}^2.
$$
The lemma is proved.

\quad

Note by the way, that for the identity homeomorphism $h_0(t)\equiv t$ it is 
obvious that $\|e^{ikh_0}\|_{W_2^{1/2}(\mathbb T)}=|k|^{1/2}$. 

\quad

\textbf{Proof of Theorem 2.} Assume that $h$ is a self-homeomorphism of 
$\mathbb T$ such that $(F\ast\lambda)\circ h\in W_2^{1/2}(\mathbb T)$ for 
every measure $\lambda\in P(\mathbb T)$. Let $M(\mathbb T)$ be the Banach 
space of all real measures $\mu$ on $\mathbb T$ with the usual norm 
$\|\mu\|_{M(\mathbb T)}$ equal to the variation of $\mu$. Each $\mu\in 
M(\mathbb T)$ is a linear combination of two probability measures, so we see 
that $(F\ast\mu)\circ h\in W_2^{1/2}(\mathbb T)$ for all $\mu\in M(\mathbb 
T)$. 

Note that $W_2^{1/2}(\mathbb T)$ is a Banach space with respect to the norm 
$$
\|f\|^\circ_{W_2^{1/2}(\mathbb 
T)}=\bigg(\sum_{k\in\mathbb Z}|\widehat{f}(k)|^2(|k|+1)\bigg)^{1/2}.
$$ 
We have $W_2^{1/2}(\mathbb T)\subseteq L^2(\mathbb T)$ and 
$\|\cdot\|_{L^2(\mathbb T)}\leq \|\cdot\|^\circ_{W_2^{1/2}(\mathbb T)}$, 
where $L^2(\mathbb T)$ is the space of all square integrable functions with 
$$
\|f\|_{L^2(\mathbb T)}=\bigg(\frac{1}{2\pi}\int_{\mathbb T}|f(t)|^2 dt\bigg)^{1/2}.
$$ 
 
Consider the (linear) operator $Q : M(\mathbb T)\rightarrow W_2^{1/2}(\mathbb 
T)$ defined by 
$$
Q\mu=(F\ast\mu)\circ h.
$$ 
We claim that $Q$ is a bounded operator. To verify this consider a sequence 
$\mu_n\in M(\mathbb T),\, n=1, 2,\ldots,$ such that $\mu_n\rightarrow\mu$ in 
$M(\mathbb T)$ and $Q\mu_n\rightarrow g$ in $W_2^{1/2}(\mathbb T)$. 
Obviously, the sequence $F\ast\mu_n, \, n=1, 2, \ldots,$ converges uniformly 
to $F\ast\mu$. So the sequence $(F\ast\mu_n)\circ h, \, n=1, 2, \ldots,$ 
converges uniformly to $(F\ast\mu)\circ h$, whence $Q\mu_n\rightarrow Q\mu$ 
in $L^2(\mathbb T)$. At the same time, $Q\mu_n\rightarrow g$ in $L^2(\mathbb 
T)$. So $Q\mu=g$. By the closed graph theorem $Q$ is bounded. 

Thus, for every $\mu\in M(\mathbb T)$ we have 
$$
\|(F\ast\mu)\circ h\|_{W_2^{1/2}(\mathbb T)}\leq\|(F\ast\mu)\circ h\|^\circ_{W_2^{1/2}(\mathbb T)}
\leq c\|\mu\|_{M(\mathbb T)},
$$
where $c>0$ is independent of $\mu$. Applying this estimate to 
$\mu=\delta_{-\theta}$, where $\theta$ is a point in $\mathbb T$ and 
$\delta_{-\theta}$ is the unit mass at $-\theta$, and taking into account 
that $F_\theta=F\ast\delta_{-\theta}$, we obtain 
$$
\|F_\theta\circ h\|_{W_2^{1/2}(\mathbb T)}\leq c \quad \textrm{for all} \quad \theta\in \mathbb T.
\eqno(11)
$$

Note now that for each $k\in\mathbb Z$ one can find a positive integer $m(k)$ 
so that 
$$
\sum_{|\nu|\leq m(k)} |\widehat{e^{ikh}}(\nu)|^2|\nu|\geq |k|/2.
$$
Indeed, if $e^{ikh}\in W_2^{1/2}(\mathbb T)$, the existence of  $m(k)$ 
follows from Lemma 3, while if $e^{ikh}\not\in W_2^{1/2}(\mathbb T)$, the 
existence of  $m(k)$ is obvious. 

Let $N$ be an arbitrary positive integer. We define $M(N)$ by 
$$
M(N)=\max_{|k|\leq N} m(k).
$$ 
Then, for all $k\in\mathbb Z$ with $|k|\leq N$ we have 
$$
\sum_{|\nu|\leq M(N)} |\widehat{e^{ikh}}(\nu)|^2|\nu|\geq |k|/2.
\eqno(12)
$$ 
At the same time, Lemma 2 implies 
$$
\frac{1}{2\pi}\int_\mathbb T |\widehat{F_\theta\circ h}(\nu)|^2
d\theta\geq\sum_{|k|\leq N}|\widehat{F}(k)|^2|\widehat{e^{ik h}}(\nu)|^2,
$$
whence, by multiplying by $|\nu|$ and summing over $|\nu|\leq M(N)$, we 
obtain 
$$
\frac{1}{2\pi}\int_\mathbb T \bigg(\sum_{|\nu|\leq M(N)}|\widehat{F_\theta\circ h}(\nu)|^2|\nu|\bigg)
d\theta\geq\sum_{|k|\leq N}\bigg(|\widehat{F}(k)|^2\sum_{|\nu|\leq M(N)}|\widehat{e^{ik h}}(\nu)|^2|\nu|\bigg).
$$
Taking (11) and (12) into account, we see that
$$
c^2\geq\sum_{|k|\leq N}|\widehat{F}(k)|^2(|k|/2).
$$
Since $N$ was chosen arbitrarily, this yields $F\in W_2^{1/2}(\mathbb T)$, 
which contradicts the assumption of the theorem. The theorem is proved. 

\quad

\textbf{4. The class $\mathrm{Lip}_{1/2}(\mathbb T)$.} Here we obtain the 
main result of the paper. Recall that if $\alpha>1/2$, then 
$\mathrm{Lip}_\alpha(\mathbb T)\subseteq W_2^{1/2}(\mathbb T)$ (see (3)). On 
the other hand, the function 
$$
F(t)=\sum_{n\geq 0}2^{-n/2}e^{i2^nt}
$$
is in $\mathrm{Lip}_{1/2}(\mathbb T)$ (see, e.g., [1, Ch. XI, Sec. 6]) but is 
not in $W_2^{1/2}(\mathbb T)$. Obviously, for the corresponding family of 
translations $S_F$ we have $S_F\subseteq \mathrm{Lip}_{1/2}(\mathbb T)$ (the 
same is true for $K_F$). Thus, an immediate consequence of Theorem 1 (as well 
as that of Theorem 2) is the following Theorem 3. 

\quad

\textbf{Theorem 3.} \emph{There does not exist a self-homeomorphism $h$ of 
$\mathbb T$ such that $f\circ h\in W_2^{1/2}(\mathbb T)$ for every $f\in 
\mathrm{Lip}_{1/2}(\mathbb T)$.} 

\quad

\quad

\textbf{5. Remarks and open problems.} 1. Given two real-valued functions $u$ 
and $v$ in $\mathrm{Lip}_{1/2}(\mathbb T)$, does there exist a change of 
variable $h$ such that $u\circ h$ and $v\circ h$ are in $W_2^{1/2}(\mathbb 
T)$\,? (This question was already mentioned in the introduction.) 

2. Given a real-valued $F\in C(\mathbb T)$ and a $\theta\in \mathbb T, 
\,\theta\neq 0$, does there exist an $h$ such that $F\circ h$ and 
$F_\theta\circ h$ are in $W_2^{1/2}(\mathbb T)$\,? What if $F\in 
\mathrm{Lip}_{1/2}(\mathbb T)$\,? 

3. It is unclear if one can replace $C^+(\mathbb T)$ with $C(\mathbb T)$ in 
Theorem 1. 

4. For $s>0$ the Sobolev space $W_2^s(\mathbb T)$ is defined as the space of 
square integrable functions $f$ on $\mathbb T$ with $\sum_{k\in\mathbb Z} 
|\widehat{f}(k)|^2|k|^{2s}<\infty$. It was shown in [5, Corollary 3] that if 
$K$ is a compact set in $C(\mathbb T)$, then there exists a 
self-homeomorphism $h$ of $\mathbb T$ such that $f\circ 
h\in\bigcap_{s<1/2}W_2^s(\mathbb T)$ for every $f\in K$. 

5. There exists a real-valued $f\in C(\mathbb T)$ such that $f\circ 
h\notin\bigcup_{s>1/2}W_2^s(\mathbb T)$ whenever $h$ is a self-homeomorphism 
of $\mathbb T$. This is a simple consequence of the inclusion 
$\bigcup_{s>1/2}W_2^s\cap C(\mathbb T)\subseteq A(\mathbb T)$, where 
$A(\mathbb T)$ is the Wiener algebra of absolutely convergent Fourier series, 
and the result obtained in [10] (see also [11, Theorem 3.2]): there exists a 
real-valued $f\in C(\mathbb T)$ such that $f\circ h\not\in A(\mathbb T)$ for 
every self-homeomorphism $h$.

6. We now provide a very short proof of the refined version of the 
Bohr--P\'al theorem. Implicitly this proof is contained in the proof of 
Theorem 4 in [5]. The refinement is achieved through a slight modification of 
the original argument. Suppose that $f\in C(\mathbb T)$ is real-valued. 
Without loss of generality we assume that $f(t)>0$ for all $t\in\mathbb T$. 
Consider the curve $\gamma$ in the complex plane $\mathbb C$ given by 
$\gamma(t)=f(t)e^{it}, \,t\in [0, 2\pi]$. This is a closed continuous curve 
without self-intersections. By $\Omega$ we denote the interior domain bounded 
by $\gamma$. Consider a conformal mapping $G$ of the unit disc 
$D=\{z\in\mathbb C : |z|<1\}$ onto $\Omega$. As is well known, $G$ extends to 
a homeomorphism of the closure $\overline{D}$ of $D$ onto the closure 
$\overline{\Omega}$ of $\Omega$ and, being thus extended, homeomorphically 
maps the circle $\partial D=\{z\in\mathbb C : |z|=1\}$ onto the boundary 
$\partial \Omega$ of $\Omega$. We retain the notation $G$ for the extension 
and consider the function $g$ defined by $g(t)=G(e^{it})$. Clearly, there is 
a self-homeomorphism $h$ of the segment $[0, 2\pi]$ such that 
$g(t)=\gamma(h(t)), \,t\in[0, 2\pi]$. Note that $\pi\sum_{n\geq 
0}|\widehat{g}(n)|^2n$ is the area of $\Omega$. So $\gamma\circ h\in 
W_2^{1/2}(\mathbb T)$. At the same time, we have $f\circ h=|\gamma\circ h|$. 
It remains to observe that for an arbitrary function $F$ the condition $F\in 
W_2^{1/2}(\mathbb T)$ implies $|F|\in W_2^{1/2}(\mathbb T)$, which is clear 
from the definition of $\||\cdot|\|_{W_2^{1/2}(\mathbb T)}$ (see (3)). Thus, 
$f\circ h\in W_2^{1/2}(\mathbb T)$. 

A stronger result, based on the Riemann's theorem on conformal mappings, was 
obtained in [2], see also [8].

\quad 

\begin{center}
\textbf{References}
\end{center}

\flushleft
\begin{enumerate}

\item  N. K. Bary, \emph{A treatise on trigonometric series},
    Vols. I, II, Pergamon Press, Oxford 1964.

\item  W. Jurkat,  D. Waterman, \emph{Conjugate functions and
    the Bohr--P\'al theorem}, Complex Variables 12 (1989),
    67--70.

\item J.-P. Kahane, \emph{Quatre le\c cons sur les hom\'eomorphismes du 
    circle et les s\'eries de Fourier}, in: Topics in Modern Harmonic 
    Analysis, Vol. II, Ist. Naz. Alta Mat. Francesco Severi, Roma, 1983, 
    955--990. 

\item  J.-P. Kahane,  Y. Katznelson, \emph{S\'eries de Fourier des 
    fonctions born\'ee}, Studies in Pure Math., in Memory of Paul Tur\'an, 
    Budapest, 1983, pp.~395--410. (Preprint, Orsay, 1978.) 

\item V. V. Lebedev, \emph{Change of variable and the rapidity of decrease 
    of Fourier coefficients}, Matematicheski\v{\i} Sbornik, 181:8 (1990), 
    1099--1113 (in Russian). English transl.: Mathematics of the 
    USSR-Sbornik, 70:2 (1991), 541--555. English transl. corrected by the 
    author is available at: https://arxiv.org/abs/1508.06673 
    
\item V. V. Lebedev, \emph{Torus homeomorphisms, Fourier coefficients, and 
    integral smoothness}, Russian Mathematics (Iz. VUZ), 36:12 (1992), 
    36--41. 

\item V. Lebedev, \emph{The Bohr--P\'al theorem and the Sobolev space 
    $W_2^{1/2}$}, Studia Mathematica, 231:1 (2015), 73--81. 

\item  V. V. Lebedev, \emph{A short and simple proof of the 
    Jurkat--Waterman theorem on conjugate functions}, Functional Analysis 
    and Its Applications, 51:2 (2017), 148--151. 

\item V. Lebedev,  A. Olevskii, \emph{Homeomorphic changes of variable and 
    Fourier multipliers}, Journal of Mathematical Analysis and 
    Applications, \textbf {481}:2 (2020) 123502, 1--11. 

\item A. M. Olevski\v{\i}, \emph{Change of variable and absolute 
    convergence of Fourier series}, Soviet Math. Dokl., 23 (1981), 76--79. 

\item A. M. Olevski\v{\i}, \emph{Modifications of functions and Fourier 
    series}, Russian Math. Surveys, 40 (1985), 181--224. 

\item A. A. Saakjan, \emph{Integral moduli of smoothness and the Fourier 
    coefficients of the composition of functions}, Mathematics of the 
    USSR-Sbornik, 38 (1981), 549--561.

\end{enumerate}

\quad

\nopagebreak 

\noindent {\small School of Applied Mathematics, National Research University 
Higher School of Economics (HSE University), 34 Tallinskaya St., Moscow, 123458, Russia\\
e-mail: \emph{lebedevhome@gmail.com} }

\end{document}